\documentclass[%
 aip,cha,
 amsmath,amssymb,
 reprint,%
]{revtex4-1}

\usepackage{latexsym,graphicx,shortvrb,color,graphics,latexsym,amsfonts,amsmath,epsfig}
\usepackage{xcolor}

\def\beq{\begin{equation}}
\def\eeq{\end{equation}}

\def\Dx{\Delta x}
\def\Dy{\Delta y}
\def\Ds{\Delta s}

\def\LR{\mathbb R}

\def\xp{x^{+}}

\def\norm#1{\Vert#1\Vert}
\def\t{\theta}
\def\a{\alpha}
\def\bw{\bar{w}}

\begin{document}

\title[Chaotic behavior of the AFS algorithm]{On the chaotic behavior of the Primal--Dual Affine--Scaling Algorithm for Linear Optimization}

\author{H. Bruin}
\email{henk.bruin@univie.ac.at}
 \affiliation{{Faculty of Mathematics,
University of Vienna,
    Oskar Morgensternplatz 1\\
A-1090 Vienna, Austria }}
\author{R. Fokkink}%
 \email{r.j.fokkink@tudelft.nl}
\affiliation{
Delft University, Faculty of Electrical Engineering,
Mathematics and Computer Science,\\ P.O.Box 5031, 2600 GA
Delft, Netherlands}%
\author{G. Gu}
\email{ggu@nju.edu.cn}
\affiliation{
Department of Mathematics,
Nanjing University,
Nanjing 210093,
China
}%
\author{C. Roos}
 \email{c.roos@tudelft.nl}
\affiliation{
Delft University, Faculty of Electrical Engineering,
Mathematics and Computer Science,\\ P.O.Box 5031, 2600 GA
Delft, Netherlands}%

\date{\today}

\begin{abstract}
\noindent \normalsize
We study a one-parameter family of quadratic maps, which
serves as a template for interior point methods.
It is known that such methods
can exhibit chaotic behavior, but this has been verified
only for particular linear optimization problems.
Our results indicate that this chaotic behavior is generic.
\end{abstract}


\maketitle


We study a one-parameter family of quadratic maps on a projective simplex,
which has been derived from an interior point method, known as
the primal-dual Affine Scaling method~\cite{JansenRoosTerlaky}. This
particular method neatly handles both the primal
and the dual variables in one step, enabling us to derive a one-parameter
family, independently of the underlying linear optimization problem.
We study the bifurcations of this one-parameter family and
find that they are almost identical to those that
have previously been found by Castillo and Barnes~\cite{Barnes} for a specific linear optimization problem,
using another interior point method.
This indicates, experimentally and non-rigorously, that the route to chaos
in our one-parameter family is typical for general interior point methods.

\section{Introduction}

In linear optimization one wants to compute the maximum value of
a linear objective function under linear inequality constraints. There exist
many algorithms that solve LO problems by iteration. The classical
algorithm is the simplex method, which produces an exact
solution.
It runs through the extremal points of the convex set
that satisfies the constraints (the feasible set), improving the
value of objective function in each step, halting at an extremal point
that produces the maximum value.
The simplex method runs from one boundary point of the feasible set to the next.
Interior point methods run
through the interior of the feasible set.
Historically, the first such method is the {\sl affine
scaling algorithm} (AFS) method of Dikin, which remained unnoticed  until 1985. The
work of Karmarkar~\cite{int:Karmarkar2} sparked a
large amount of research in polynomial--time methods for LO,
and gave rise to many new and efficient interior point methods
(IPMs) for LO. For a survey of this development we refer to the
books of Wright~\cite{int:SWright9}, Ye~\cite{BookYe},
Vanderbei~\cite{Vanderbei} and Roos \emph{et
al.}~\cite{RoosTerlakyVial2}.
An IPM starts from an arbitrary initial point $x_0$ in the interior
and constructs a sequence $x_n$
that converges to a maximum $x^*$.
An IPM is a dynamical system
which solves the LO problem, provided that the $\omega$-limit of~$x_0$
consists of maxima of the objective function~$f$.
If there is only one such maximum,
then the LO problem is called \textit{non-degenerate}.
In this case, the IPM solves the LO problem provided that it converges
to the maximum.

Any LO problem can be converted to a dual problem in which
one needs to find the minimum value of a dual linear objective
function under dual constraints. If $f$ is the objective function
of the primal problem and if $g$ is the objective function
of the dual problem, then
$f(x)\leq g(y)$ for all feasible $x$ and
all $y$.
To solve an LO problem it therefore suffices to close the
duality gap and find $x^*,y^*$ such that $f(x^*)=g(y^*)$.
According to the minimax theorem, such $x^*,y^*$ exist
and apart from a primal sequence $x_n$
most IPM's also produces a dual sequence $y_n$,
halting as soon as the duality gap $g(y_n)-f(x_n)$
reaches a value which is below the desired accuracy threshold.
A primal problem that is non-degenerate
may have a degenerate dual problem, in which case the orbit
$y_n$ may have a non-trivial $\omega$-limit set.
We will study such an LO problem at the end of this paper and
find that the dual dynamical system
contains a hyperbolic attractor.

The simplex method runs from one extremal point to the next,
but an IPM uses a variable step size~$\alpha$. For each feasible $x$ the
algorithm produces a vector $v$ such that $x+\alpha v$
is contained in the feasible set for $0\leq \alpha\leq 1$
and such that $x+v$ is in the boundary of the feasible set.
The step size $\alpha$ is fixed during iteration and
is chosen $<1$ so that the orbit $x_n$ is contained in
the interior of the feasible set. An IPM therefore is a one-parameter
family of dynamical systems. It is well known that
an IPM may not converge if $\alpha$ is too large.
One of the best studied algorithms is the Affine Scaling method
(AFS) which was proposed by Dikin and which has been further
developed by Vanderbei {\em et al.} \cite{Vanderbeietal} It
is known~\cite{TsuchiyaMuramatsu} that AFS converges if
$\alpha\leq 2/3$ and that it need not converge if $\alpha>2/3$,
see~\cite{twothirds}. It is also known that AFS behaves
chaotically in the dual variables for $\alpha>2/3$,
as has been found by Castillo and Barnes~\cite{Barnes} and Mascarenhas \cite{Mascarenhas}.

\subsection{Outline of our paper}

The previous studies of chaotic behavior
in interior point methods were carried out for specific
problems: one considers an LO problem, applies the algorithm
and analyzes the resulting dynamical system.
In this paper, we take a different approach.
We consider the primal-dual AFS
that was proposed by Jansen {\em et al.}~\cite{JansenRoosTerlaky}.
It has the nice property that
it be presented in a such a form that
its low order terms do not depend on the original LO problem.
By ignoring the
higher order terms, we obtain a one-parameter family of dynamical systems,
which we call the \textit{Dikin process},
that is the same for all LO problems.
Of course, the Dikin process
is not an IPM anymore.
However, the bifurcations that we establish for the Dikin process
are the same as
the bifurcations that have previously been found by Castillo and Barnes for
their specific LO problem. This indicates, experimentally
and non-rigorously, that the chaotic behaviour of the Dikin process
represents that of general interior point methods.

Our paper is organized as follows. We first recall the primal-dual
AFS method for solving LO problems. We then derive the one-parameter
family of dynamical systems, and analyze it for increasing values of
a parameter $\t$. We show that the system behaves chaotically as $\t$
increases beyond $2/3$. We supplement this analysis experimentally
by Feigenbaum
diagrams. In the final section, we compare our results to an
IPM that arises from a specific LO problem.

\subsection{Notation}
We reserve the symbol $e \in \LR^n$ for the vector of all
ones. For a vector $x$, the capital $X$ denotes the diagonal
matrix with the entries of $x$ on the diagonal. Furthermore, if
$f:\LR\rightarrow\LR$ is a function and $x \in \LR^n$, then we
denote by $f(x)$ the vector $(f(x_1),\dots,f(x_n))$. If $s$ is
another vector, then $xs$ will denote the coordinatewise product of
$x$ and $s$ and $x/s$ will denote the coordinatewise quotient of $x$
and $s$. In other words, $xs=Xs$ and $x/s=S^{-1}x$.  Finally, $\| .
\|$ denotes the $l_2-$norm.

\section{A recap of the Primal-dual affine scaling method}
In linear optimization, the notion of affine scaling has been introduced by Dikin~\cite{int:Dikin1}
as a tool for solving the (primal) problem in
standard format
\[
 (P)\;\;\;\qquad \min\{c^T x\,:\, A x=b x\ge0\}.
\]
The underlying idea is to replace the nonnegativity constraints
$x\ge 0$ by the ellipsoidal constraint
\beq\label{ellips} \norm{\bar{X}^{-1}({\bar x}-x)}\leq 1,\eeq
where $\bar{x}$ denotes some given interior feasible point, and
$\bar{X}$ the diagonal matrix corresponding to $\bar{x}$. The
resulting subproblem is easily solved and renders a new interior
feasible point with a better objective value. Dikin showed, under
the assumption of primal nondegeneracy, that this process converges
to an optimal solution of $(P)$.

Every known method for solving $(P)$ essentially also solves the
dual problem
\[
 (D)\;\;\;\qquad \max\{b^T y\,:\, A^T y + s=c s\ge0\}
\]
by closing the duality gap between $c^Tx$ and $b^Ty$, which equals
$x^Ts$.  Our basic assumption is that a
primal-dual pair $(x,s)$ of feasible solutions exists and that $A$
is an $m\times n$ matrix of rank $m$ for $m<n$. A pair
of feasible vectors $x^*,s^*$ solves $(P)$ and $(D)$ if and only if they are
orthogonal. Since $x^*\geq 0$ and $s^*\geq 0$, this means that the
coordinatewise product $x^*s^*$ is equal to the all-zero vector.
The primal-dual AFS method that we consider in this paper has been
proposed by Jansen {\em et al.}~\cite{JansenRoosTerlaky}.
In primal-dual AFS, Dikin's ellipsoidal constraint
(\ref{ellips}) is
replaced by a constraint that includes both
the primal and the dual variable:

\beq\label{ellips2} \norm{\bar{X}^{-1}({\bar x}-x)+\bar{S}^{-1}({\bar s}-s)}\leq 1,\eeq
where $\bar {S}$ denotes the diagonal matrix
corresponding to the slack vector $\bar s$. In this notation, $(\bar
x,\bar s)$ is the original pair of primal vector and slack vector
and $(x,s)$ is an updated pair. The differences $\Delta x=x-\bar x$
and $\Delta s=s-\bar s$ are called the primal-dual AFS directions.

For non-negative $x,s$ let $v=(xs)^{1/2}$ be the coordinatewise square root of the coordinatewise
product and let $v^k$ be the coordinatewise power of $v$.
Jansen {\em et al.}\ have shown that the directions $\Delta x$ and $\Delta s$
can be derived from the vector
\[p_v=-\frac{v^{3}}{\norm v^2}\]
by first projecting $p_v$ onto the null space (for $\Delta x$) and the row space (for $\Delta s$)
of $AD$, with $D=XS^{-1}$, and then rescaling the result by a coordinatewise product.
More specifically, if $d=(x/s)^{1/2}$ then
\[\Delta x=dP_{AD}(p_v),\ \Delta s= d^{-1}Q_{AD}(p_v)\]
where $P_{AD}$ and $Q_{AD}$ denote the orthogonal projections onto
the null space of $AD$ and the row space of $AD$, respectively.
These projections recombine in the Dikin ellipsoid to
\beq\label{sum} X^{-1}\Delta x + S^{-1}\Delta
s=-\frac{v^2}{\norm{v^2}}. \eeq This gives the primal-dual AFS directions but
not the size of the step, which is controlled by an additional
parameter $\alpha$. It is known that
the iterative process $x\mapsto
x+\alpha\Delta x, s\mapsto s+\alpha\Delta s$  converges to a
solution if
$\alpha<1/\left(15\sqrt n\right)$.\cite{JansenRoosTerlaky} This is of course a significant restriction on
the step size, and primal-dual AFS is not often used in practice.

\section{Derivation of the Dikin process}

Starting with a primal-dual feasible pair $(x,s)$, the next iterated
pair is given by
\[\xp = x+\a\Dx  = s+\a\Ds,\]
and hence we have
\[\xp = xs+\a\left({x\Ds+s\Dx}\right) + \a^2\Dx\Ds.\]
The vectors $\Dx$ and $\Ds$ are orthogonal. If
the AFS iterations are close to a solution,
then $\Dx$ and $\Ds$ will be relatively small, and the product
$\Dx\Ds$ will be negligible. If we ignore the
quadratic term, i.e., if we assume that the coordinatewise product
$\Dx\Ds$ is equal to zero,
then the reduction of~$xs$ is proportional to~$x\Ds+s\Dx$, which can
be rewritten to
\[
xs\left(x^{-1}\Delta x+ s^{-1}\Delta s\right).
\]
Observe that $x^{-1}\Delta x+ s^{-1}\Delta s$ is equal to the
left-hand side in $(\ref{sum})$.
So if we ignore the quadratic term, and if we use equation~$(\ref{sum})$,
then we find that
\[\xp  = xs+\a\left({x\Ds+s\Dx}\right) = xs -\a \frac{x^2s^2}{\norm{xs}} = xs{e -\a \frac{xs}{\norm{xs}}}.\]
Recall that $e$ denotes the all-one vector, so we may also write this as
\[
\xp =xs{e -\a \frac{xs}{\norm{xs}}}.
\]
Now we have arrived at an iterative process for the
product vector $xs$. Since we require $x\ge 0$ and $s\ge 0$, we need
to require $xs\ge 0$ in the iterative process, and the maximal step size is equal to
\[\a_{\max} = \frac{\norm{xs}}{\max{xs}}\]
Defining
\[\t = \frac{\a}{\a_{\max}} = \a\frac{\max{xs}}{\norm{xs}}\]
and writing $w=xs$ we get
\[w^+ =  w{e -\t \frac{w}{\max{w}}}\quad \t\in [0,1].\]
This iterative process depends on a parameter $\t$ which is related
to the original step size by $\alpha=\t\a_{max}$.
If $w$ has coordinates that are
approximately equal (in optimization one says that $w$ is `close to the central line'),
then $\a_{max}\approx \sqrt n$. In general,
$1\leq \a_{max}\leq \sqrt n$.

We make one further reduction. If $u=\lambda w$ for a scalar
$\lambda$ then $u^+=\lambda w^+$, so the iterative process preserves
projective equivalence. We may therefore reduce our system up to
projective equivalence by scaling vectors so that their maximum
coordinate is equal to one. If we consider vectors up to projective
equivalence, then we obtain our \textit{Dikin process}:
\beq\label{eq:scaled process}
{\bw}^{k} = w^k{e -\t w^k}\quad w^{k+1} =
\frac{{\bw}^{k}}{\max{\bw}^{k}}\quad k=0, 1, \ldots
\eeq
The Dikin process involves two
steps: multiplication and scaling. To describe the process more succinctly
we use the map $f_\t(x)=x(1-\t x)$. The
Dikin process is then given by:
\begin{equation}
w^{k+1}=f_\t(w^k)/\max\{f_\t(w^k)\}.\end{equation}
Note that $f_{\t}$ is a higher-dimensional analog of the logistic map on the unit interval.
For each coordinate we apply the same quadratic map,
and the only interaction between the coordinates is induced by the scaling.
The Dikin process does not solve the original LO problem. Its significance
derives from the fact that it does not depend on the LO problem and that its
bifurcations can be analyzed in a standard way.

\section{Bifurcation analysis}

We analyze the Dikin process $f_\t$,
for increasing values of $\t$.
We suppress the subscript $\theta$ in $f_\t$ and simply write the
Dikin process as
\begin{equation}w^{k+1}=f(w^k)/\max\{f(w^k)\}.\end{equation}
Note that $f$ has a global maximum $f(1/2\t)=1/4\t$
and that we scale $f$ such that all coordinates take values $\leq 1$.

\subsection{$\t\leq 2/3$: the process converges to $e$}
If
$\t\leq 1/2$ then the global maximum of $f$ is $\geq 1$,
which is outside the domain of our coordinates. The value
of each coordinate increases during iteration.
By monotonicity the limit of $w^k$ exists and it is a
fixed point under iteration. The only fixed point is $e$ and
therefore $w^k$ converges to the all-one vector~$e$ if $\t\leq 1/2$.

We now argue that $e$ remains the global attractor if $\t\leq 2/3$.
If $\t> 1/2$ then $f$ is unimodal
and point symmetric with respect to its maximum $1/2\t$:
\beq\label{pointsym} f\left(\frac 1{2\t}+z\right)=f\left(\frac
1{2\t}-z\right).\eeq
Under iteration of $f$ all orbits eventually end up in
the interval $[1-1/\t,1]$. In particular,
for every initial $w^0$ it eventually holds that
$\min{w^k}\geq 1/\t-1$.
We need to show that $\min w^k$
in fact converges to $1$ if $1/2\leq \t\leq 2/3$.
By the point symmetry in
\eqref{pointsym}, if we replace the coordinates $w_i<1/2\t$ in $w^k$
by their reflections $1/\t-w_i$, then this does not affect
$w^{k+1}$. We may therefore assume that $\min w^k\geq 1/2\t$.
Let $x=\min w^k\geq 1/2\t$. Since $f$ is
decreasing for $x\geq 1/2\t$ we have that $f(x)=\max
f(w^k)$ and $f(1)=\min f(w^k)$. Therefore the minimum coordinate of
$w^{k+1}$ is given by
\begin{equation}\label{h}
h(x)=\frac {f(1)}{f(x)}=\frac{1-\t}{x(1-\t x)}.
\end{equation}
To prove that the process converges to $e$, it now
suffices to prove that $h(x)\geq x$, because this implies that
the limit of $h^k(x)$ exists and is equal to the unique fixed point
of $h$. Now $h(x)\geq x$ can be rewritten as \beq\label{cbc} \t
x^3-x^2+1-\t\geq 0. \eeq
The derivative $(3\t x - 2)x$ of the cubic $\t x^3 -x^2 +1 -\t$ is
negative on the unit interval, by our assumption that $\t\leq 2/3$.
So the cubic has its
maximum at $0$ and its minimum at $x=1$, which is a zero of the cubic.
Hence the inequality $h(x)\geq x$ holds and we conclude that $w^k$
also converges to $e$ if $1/2\leq\t\leq 2/3$.

\subsection{$2/3<\t\leq \frac{1+\sqrt
5}4$: convergence to a point of period two.}

We will see that if $2/3<\t\leq \frac{1+\sqrt
5}4$, then the
minimum coordinate and the maximum coordinate interchange
under iteration, while all other coordinates either
converge to the minimum of the maximum.
We can thus ignore these other coordinates, and
observe that the Dikin
process on the minim and maximum coordinate is given by $(x,1)\to(1,h(x))\to (h^2(x),1)\to\cdots$,
with $h$ as in~(\ref{h}).
Observe that a fixed point of $h$ produces a point of period two for this
process.

If $\t>2/3$ then $h$ has a unique fixed point
$r\in (0,1)$, which can be found by solving the cubic equation
$h(x)=x$ that we already encountered in equation~\eqref{cbc}. This
cubic is divisible by $x-1$, so we find that
$r$ satisfies the quadratic equation \beq
\label{quadratic} \t r^2+(\t-1)r+(\t-1)=0.\eeq
The positive solution for $r$ is equal to
\beq\label{eq:r}
r=\frac{1-\t+\sqrt{(1-\t)^2+4\t(1-\t)}}{2\t},
\eeq
which is $\leq 1$ if and only if $\t\geq 2/3$.
The cubic equation $h(x)=x$ has zeros in $1,r$ and the third zero
$s$ is negative. In particular, $h(x)>x$ on $(s,r)$ and $h(x)<x$ on
$(r,1)$ and we find that $r$ is the global attractor
of $h$ in the interval $[1/\t-1,1)$.
Note that $r$ is not a global attractor in the closed interval
$[1/\t-1,1]$ since $h(1)=1$ is a fixed point.

Since $r$ is an attractor, the two-dimensional process
$(x,1)\to(1,h(x))$ converges to an orbit of period two if $\t> 2/3$.
We note that this
particular limit behavior has also been observed by Hall and Vanderbei
\cite{twothirds} for the (primal) AFS method.

Since all coordinates eventually increase above $1-1/\t$
we may as well assume that $\min w^0\geq 1/\t-1$.
The minimum coordinate of $f(w^0)$ then has value
$f(1)=1-\theta$ and the maximum coordinate has value $\leq 1/4\t$.
Therefore, $\min w^1\geq 4\t(1-\t)$. If in fact $\min w^1\geq 1/2\t$,
then we can restrict out attention to the mimimum and the maximum
coordinate. This would be the case if $4\t(1-\t)\geq 1/2\t$, which leads
us to the cubic equation
$$8\t^2(1-\t)-1=0\Longleftrightarrow (2\t-1)(1-2\t-4\t^2)=0.$$
The two roots of the quadratic are $\frac{1\pm\sqrt 5}4$, and
so we conclude that we may indeed restrict our attention to the
minimum and the maximum coordinate if $\frac 12 < \t\leq
\frac{1+\sqrt 5}4$.

Suppose $\t\leq \frac {1+\sqrt 5}4\approx 0.809$ and consider an
initial condition $w^0 = (w_1, \dots, w_n)$ with increasing
coordinates $w_1\leq w_2 \leq \dots \leq w_n = 1$ and such that
$w_1\geq 1/2\t$.  For an intermediate
coordinate $w_i$ the process is given by $ w_i \to
\frac{w_i(1-\theta w_i)}{w_1(1-\theta w_1)}. $ Now the minimum coordinate
converges to $r$ so we may as well put $w_1=r$, in which case
we get that $w_i\to g(w_i)$ for the map
$$
g(x) = \frac{x(1-\theta x)}{r(1-\theta r)}.
$$
This map $g$ keeps track of the Dikin process $w^k$ on a fixed coordinate.
We now prove that the $\omega$-limit of $g$ is Lebesgue a.e. equal to $\{r,1\}$,
which will show that a.e. point converges to a point of period two.
This comes down to a straightforward computation, which is carried out
in the paragraph below.

First note that $g(r)=1$ and that $g(1)=h(r)=r$ so that $g$ has a fixed point
$s\in (r,1)$ as is illustrated by the graph of $g^2$ in the figure
below. We leave it to the reader to verify that $g$ has a
unique fixed point in $(r,1)$ at $s=(r +\t - 1)/r\t$.
\begin{figure}[!hbt]
\resizebox{2in}{!}{\includegraphics[width=3in]{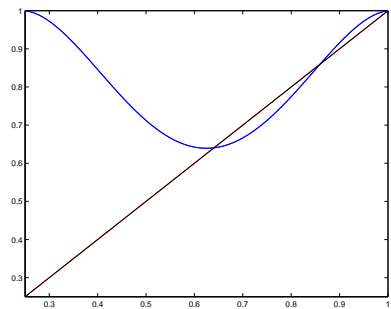}}
\caption{The function $g^2$ on the interval [1/4,1] for
$\theta=0.8$ plotted against the diagonal.}\label{ggraph}
\end{figure}
The derivative of $g$ is given by
\beq\label{derivative}
g'(x) = \frac{1-2\theta x}{r(1-\theta r)}=\frac{r(1-2\theta
x)}{1-\t},
\eeq
where we use that $1-\t=r^2-\t r^3$.
Note that $g(x)=x$ is a quadratic equation with solutions $x=0$ and $x=s$.
The
equation $g^2(x)=x$ is an equation of degree four with solutions $0,r,s,1$. It follows
that $g^2(x)\not=x$ on the two subintervals $(r,s)\cup(s,1)$ and
that $g^2(x)>x$ on the one interval while $g^2(x)<x$ on the other
interval. Using equation~(\ref{derivative}), and using that
$rs\t=r+\t-1$, we find that the derivative at
$s$ is
\begin{eqnarray*}
g'(s) &=&  \frac{r\left( 1 - 2\t s
\right)}{1-\theta}=\frac{2-r-2\t}{1-\t}.
\end{eqnarray*}
To prove that $s$ is unstable, we need to verify that
$\frac{2-r-2\t}{1-\t}<-1$, or equivalently, that $r>3-3\t$.
Substituting (\ref{eq:r}) for $r$ and simplifying equations
we end up with
$
(1-\t)+\sqrt{(1-\t)^2+4\t(1-\t)}>6\t(1-\t).
$
Taking squares to remove the root gives
$
(1-\t)^2+4\t(1-\t)>(6\t-1)^2(1-\t)^2.
$
which simplifies to
$
1+3\t>(6\t-1)^2(1-\t)
$.
Collecting all terms and dividing by $\t$ we finally arrive
at the inequality $9\t^2-12\t+4>0$, or equivalently,
$(3\t-2)^2>0$. This obviously holds if $\t>2/3$.
 It follows that $g^2(x)<x$ on
$(r,s)$ and that $g^2(x)>x$ on $(s,1)$.
This completes our computation and we conclude that
if $2/3<\t\leq \frac{1+\sqrt 5}4$, then the $\omega$-limit is Lebesgue a.e. equal
to an orbit of period two. The coordinates of these periodic points
are either equal to $r$ or $1$.

\subsection{Persistence of period two.}
For generic period doubling bifurcations in smooth
dynamical systems, the parameter curve of the periodic points of
period $2n$  is parabolic and intersects the curve of the periodic
point of period $n$ transversally. At the point of intersection, the
period $n$ point changes from stable to unstable, or vice versa.
Curiously, this scenario fails in at least the first two periodic doublings
in our Feigenbaum diagrams, in particular see Figure~\ref{bifudiag}. Our
numerical
experiments show that the period two limit cycle persists beyond
$\frac{1+\sqrt{5}}4$. It is indeed possible to prove that the period
two point persist, but the analysis gets involved. We limit
ourselves to the case that $w$ has three coordinates. Assuming that
the coordinates are ordered $x<y<1$ we can write

\begin{eqnarray*}
(x,y,1) &\mapsto&
\left( 1, \frac{ y(1-\theta y) }{ x(1-\theta x) }, \frac{ 1-\theta }{ x(1-\theta x) } \right)
\\
&\mapsto&
\left(
\frac{ x(1-\theta x) }{ 1-\theta \frac{1-\theta}{x(1-\theta x)} }\ ,\
\frac{y(1-\theta y)}{1-\theta}
\frac {1-\theta \frac{y(1-\theta y)}{x(1-\theta x)} }
{1-\theta \frac{1-\theta}{x(1-\theta x)} }\ ,\ 1 \right)
\end{eqnarray*}
so we can describe the second iterate by the function
$$
F(x,y)=\left(\frac{ x(1-\theta x) }{ 1-\theta \frac{1-\theta}{x(1-\theta x)} }\ ,\
\frac{y(1-\theta y)}{1-\theta}
\frac {1-\theta \frac{y(1-\theta y)}{x(1-\theta x)} }
{1-\theta \frac{1-\theta}{x(1-\theta x)} }\right).
$$

\begin{figure}[!hbt]
 \resizebox{3in}{!}{\includegraphics[width=2in]{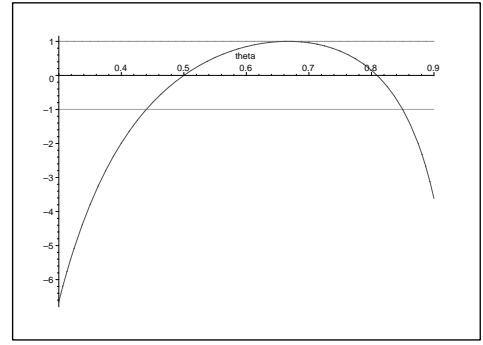}}
 \caption{\label{fig:secondeigenvalue} Value of the second
`transversal' eigenvalue $\frac{\partial F_2}{\partial y}(r,r)$ as
function of $\theta$.}
\end{figure}

This function preserves the diagonal, on which we have the
two-dimensional process, which as we have seen already has a period
two global attractor for $\theta>2/3$. So, the instability has to
occur in the direction transversal to the diagonal. We can study
this stability by taking the derivative
$$
DF(x,y) =
\left[ \begin{array}{ll}
\frac{\partial F_1(x,y)}{\partial x} & \frac{\partial F_2(x,y)}{\partial x} \\
0 &  \frac{\partial F_2(x,y)}{\partial y}
\end{array}\right]
$$
where the partial derivatives on the diagonal of the matrix $\frac{\partial F^1(x,y)}{\partial x}, \frac{\partial F_2(x,y)}{\partial y}$
are equal to
$$
\frac{(x-3\theta x^2-2\theta+2\theta^2+2\theta^2 x^3+4\theta^2x-4\theta^3x)x(1-\theta x)}
{(-x+\theta x^2+\theta-\theta^2)^2}
$$
and
$$
\frac{x-\theta x^2-2 y\theta+6 y^2\theta^2-2x^2\theta+2y\theta^2x^2-4y^3\theta^3}{(-x+\theta x^2+\theta-\theta^2) \cdot (-1+\theta) }.
$$
Maple computations show that fixed point becomes unstable at $\theta
= 0.8499377796$, when the eigenvalue $\frac{\partial
F_2(r,r)}{\partial y}$ becomes equal to $-1$. At this value of
$\theta$ we expect $(r, r,1)$ to become unstable, splitting off a
stable period $4$ point in a period doubling bifurcation, which is
confirmed by the Feigenbaum diagrams below.
In our computational results for real LO problems, we find
that the limit two cycle persists slightly beyond the threshold
of $\theta
= 0.8499377796$.

\subsection{$\t>\frac{1+\sqrt 5}{4}$: comparison to the logistic family.}

It is hard to extend the bifurcation analysis for $\t\geq
\frac{1+\sqrt 5}4$, since the degree of the algebraic equations
increases and periodic points cannot be found in closed form.
However, using the similarity between the Dikin process and
the logistic map\cite{Feigenbaum} $Q_\t:x \mapsto 4\t x(1-x)$, we can prove
that stable periodic points of higher order appear if $\t$
increases beyond~$\frac{1+\sqrt 5}{4}$. In particular,
we shall now show that if the critical point
$c = \frac12$ is $m$-periodic under $Q_\t$ then
the Dikin process
has a locally stable $m$-periodic orbit, provided the number of coordinates $n \ge m$.

Assume that the first $m$ coordinates $w_i^k$ of the vector $w^k$
are equal to $Q_\t^i(c)/\t$ (so the coordinates are not put in increasing order here). In particular, $w_m^k = c/\t = 1/2\t$,
and $f_\t(w_m^k) = 1/4\t = \max\{f_\t(x)\}$.
Then
$w_i^{k+1} = 4\t f_\t(w_i^k) = 4\t w_i^k(1-\t w_i^k)$.
The linear scaling $h(x) = \t x$ conjugates this to $Q_\t$, since
$h^{-1} \circ Q_\t \circ h(x) = 4\t x(1-\t x)$.
Since the the critical point of $Q_\t$ is periodic by our assumption,
the critical point of $f_\t$ is periodic too: $w_i^{k+1} = w_{(i \bmod m)+1}^k$ for $i = 1, \dots, m$,
and $w^{k+1}_{m-1} = 1/2\t$. In particular, the scaling remains
the same for all iterates.

This periodic orbit attracts the coordinates $w_i$ for $m < i \leq n$
and Lebesgue-a.e.\ initial choice of $w_i$.
Let us now verify that the orbit is also stable under small changes in the coordinates $w_i$ for $1 \leq i \leq m$.
Renaming these $w_i$ to $y_i$, $i=1, \dots, m$, where $y_{m-1} = 1/2\t$, $y_m = 1$, $y_1 = f(1)/f(y_{m-1})$  and $y_{i+1} = f(y_i)/f(y_{m-1})$
for $1 \leq i < m$, we can describe them by the map
\begin{equation}\label{Fmap}
F(y_1, \dots, y_{m-1}, 1) = \left( \frac{f(1)}{f(y_{m-1})},
\dots , \frac{f(y_{m-2})}{f(y_{m-1})}, 1
\right)
\end{equation}
The final coordinate is redundant, so $DF$ is an $(m-1)\times (m-1)$ matrix.
Recall that $f'(x) = 1-2\t x$. Therefore
$DF(y)$ is equal to
\[
\left( \begin{array}{ccccr}
0 & \dots & 0 &\ -(1-2\t y_{m-1})\frac{f(1)}{f(y_{m-1})^2} \\[1mm]
\frac{1-2 \t y_1}{f(y_{m-1})} & \dots & 0 &\ -(1-2\t y_{m-1})\frac{f(y_1)}{f(y_{m-1})^2} \\[1mm]
0  &  \dots & 0 &\ -(1-2\t y_{m-1})\frac{f(y_2)}{f(y_{m-1})^2} \\[1mm]
\vdots &  & \vdots & \vdots&\  \\[1mm]
0 & \dots & \frac{1-2 \t y_{m-2}}{f(y_{m-1})}  &\ -(1-2\t y_{m-1})\frac{f(y_{m-2})}{f(y_{m-1})^2}
\end{array} \right)
\]
and since $y_{m-1} = 1/2\t$, the right-most column is zero.
Therefore all eigenvalues are zero, and $DF^n$ is a contraction.
We conclude that the structure of the Feigenbaum map of the logistic
family must be present within the Feigenbaum diagrams of the Dikin process.
However, we made no estimate on the basin of attraction of the periodic points,
and our numerical results indicate that these basins are small.

\subsection{The process converges to a periodic point for $\t$ near
$1$.}\label{locstabletheta1}

Surprisingly, it is possible to determine the limit of $w^k$ for
$\t$ arbitrarily close to $1$. To conclude our bifurcation analysis,
we show that  for $\t$ close to $1$ the Dikin process
has a locally stable point of period $n$, i.e.,
the period is equal to the dimension.

Let $y$ be any point with maximal coordinate $1$ and all other
coordinates $\leq \frac 1{2\t}$. As before, we assume $\min y \geq
1/\t - 1$ and this implies that $f(1)$ is the minimal coordinate of
$w^1$. We arrange the coordinates of $y$ in
non-decreasing order. Then $f(y_{m-1})$ is the largest coordinate
among all the $f(y_k)$, so we scale by this number and we arrange
the coordinates of $w^1$ in non-decreasing order.
The dynamic process can then be described by the map $
F(y_1, \dots, y_{m-1}, 1) = \left( \frac{f(1)}{f(y_{m-1})},
 \dots , \frac{f(y_{m-2})}{f(y_{m-1})}, 1
\right)
$ of equation $(\ref{Fmap})$,
and therefore we find cyclic periodicity if
$f(y_k)/f(y_{m-1})=y_{k+1}$ and $f(1)/f(y_{m-1})=y_1$. Fix
$y_{m-1}<1/2\t$ and define a map $g(x)=f(x)/f(y_{m-1})$. Note that
$y  $ has the required cyclic periodicity if
$$
y_{m-1}=g(y_{m-2})=\cdots=g^{m-2}(y_1)=g^{m-1}(1).
$$
By the point symmetry of $f$ in \eqref{pointsym}, we may replace
$g^{m-1}(1)$ by $g^{m-1}(1/\t-1)$. If we take $y_{m-1}=1/2\t$ then a
sufficient condition for the cyclic periodic point to exist is
\beq\label{percondit}
g^{m-1}(1/\t-1)\leq 1/2\t.
\eeq
This inequality is satisfied if $\t$ is sufficiently close to $1$. Now $g$ increases
as $y_{m-1}$ decreases, so once the condition is satisfied, there
exists an $y_{m-1}$ such that $g^{m-1}(1/\t-1)=y_{m-1}$.
To compute the stability of this orbit, we cannot use anymore that the right-most column of
$DF$ vanishes, because now $y_{m-1} < 1/2\t$.
Fortunately, $DF$ is of a simple form
$$
DF(y) = \left( \begin{array}{ccccc}
0 & 0 & \dots & 0 & -c_1 \\[1mm]
\ d_1 \ & 0 &  & \vdots & -c_2 \\[1mm]
0 & \ d_2 \ &\quad  \ddots \quad &  & -c_3 \\[1mm]
\vdots &  & \ddots & & \vdots \\[1mm]
0 & \cdots & & d_{m-1} & -c_m
\end{array} \right)
$$
where $d_1>d_2>\ldots>d_{m-1}>0$ and $0 < c_1< c_2<\ldots<c_m<1$.
This follows from the fact that $y_1<y_2<\cdots<y_{m-1}\leq 1/2\t$ and that $f$ is increasing on $[y_1,y_{m-1}]\subset [0,\frac 1 {2\t}]$.
In order to estimate the eigenvalues of $DF$ we use the classical result
of~Enestr\"om-Kakeya\cite{Enestrom} that a polynomial
$p(z)=\sum_{k=0}^m a_k z^k$ with all coefficients $a_i\geq 0$
has zeros in the
annulus $\alpha\leq |z|\leq \beta$, where
\[
\alpha=\min \left\{\frac {a_k}{a_{k+1}}\right\},\ \beta=\max \left\{\frac
{a_k}{a_{k+1}}\right\}.
\]

\begin{quote} {\bf Claim:}
If $c_id_i < c_{i+1}$ for all $i \leq m-1$, then all eigenvalues of $DF$ are in the open unit
disc.
\end{quote}

Abbreviate $A = DF$ and
let $p(\lambda) = \det(\lambda I_m - A) = \sum_{k=0}^m a_k\lambda^k$
be the characteristic polynomial of $A$.
We will show by induction that the
coefficients are decreasing. More precisely $1=a_m>\cdots>a_0>0$ and
$a_0=c_1d_1\cdots d_{m-1}$.
The proof of the claim is by induction. The claim is
obvious for $m=1$. Assume that the claim is true for $m-1$. The
characteristic polynomial is equal to
$$
\lambda\det(\lambda I_{m-1}-A_{11})+(-1)^{m-1}c_1\cdot (-d_1)\cdots
(-d_{m-1}),
$$
where $A_{11}$ is the $(1,1)$-minor matrix of $A$. By the
inductive hypothesis, $\det(\lambda I_{m-1}-A)$ has decreasing
coefficients and constant coefficient $c_2\cdot d_2\cdots d_{m-1}$.
If we rewrite $a_0 = \frac{c_1}{c_2} \cdot d_1 \cdot a_1 < a_1$,
then the claim follows. We now compute
$$
\frac{c_i}{c_{i+1}} \cdot d_i
= \frac{f(y_{i-1})}{f(y_i)}\frac{1-2\t y_i}{f(y_{m-1})}
= \frac{1-2\t y_i}{1-\t y_i} < 1.
$$
which demonstrates that $c_id_i < c_{i+1}$, as claimed.
By the Enestr\"om-Kakeya Theorem, the roots of $p(\lambda)$ are all in
the open unit disc. Hence $DF^m$ is a contraction at $(y_1,\dots,y_m)$
for $\t$ sufficiently close to $1$.

Our numerical simulations suggest that the set of initial values
$w^0$ that converge to this periodic point is large and has (nearly)
full measure, as illustrated by the Feigenbaum diagrams in the next section.

\section{Feigenbaum diagrams}\label{sec:Feigenbaum}

The Dikin process is defined in aribtrary dimensions,
so in our Feigenbaum diagrams we have to project $n$-dimensional $\omega$-limits onto
one dimension. We
have chosen to simply plot one single coordinate of the $\omega$-limit set.

\begin{figure}\label{feigdim34}
{\includegraphics[height=3in]{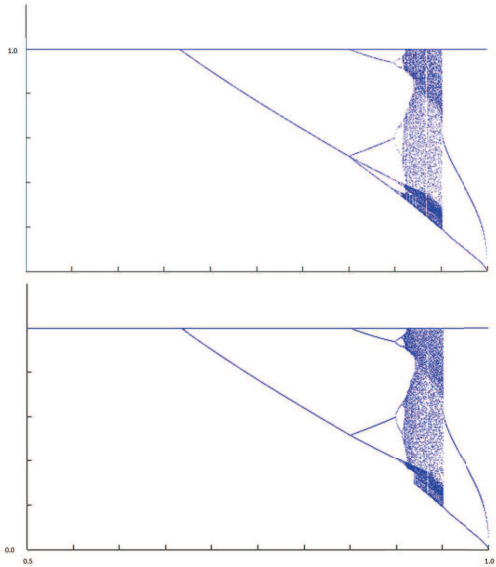}}
\caption{\label{bifudiag} Feigenbaum diagram for the process on
three coordinates. Above: $\omega$-limit of a random
coordinate, and below: $\omega$-limit of the middle coordinate.
Note that the bifurcation at 0.8499377796 produces a
period four point.}
\end{figure}

The Feigenbaum diagrams seem to exhibit the usual structure
of period doubling cascades of the logistic family
$Q_\t:x \mapsto 4\t x(1-x)$, $\t \in [0,1]$.
It is well-known that for $Q_\t$, between two period doubling bifurcations,
there is a parameter where the critical point $c = \frac12$ is periodic.
We proved above that this periodic point should then also
appear as a stable periodic point in the Dikin process, provided the dimension exceeds the period.
Since our examples have small dimension,
we do not see much of the period doubling cascade of the logistic family.

To illustrate the consequence of the choice of the projection,
compare the Feigenbaum diagrams in Figure~3. In the
top figure we plot the $\omega$-limit set of a random coordinate. Below
we choose the middle coordinate of the ordered vector. We see
that the process bifurcates at $\t=2/3$, when a point of order two
appears, and then at $\t=0.849...$, when a point of order four
appears. In the top figure, the diagram splits into five lines at
$\t=0.849...$, in the figure below it splits into four lines. The
reason for this is that the point of order four is of the type
$$
(r,s_2,1)\to(1,s_3,s_1)\to(r,1,s_2)\to(1,s_1,s_3)\to(r,s_2,1)
$$
for values $s_1,s_2$ close to $r$ and $s_3$ close to $1$. We will
plot the diagrams in the same way as the figure above, so the
reader should keep in mind that, contrary to standard Feigenbaum
diagrams, the period of a point may be smaller than the number
of lines.

The diagram indicates that $\omega$-limit set gets positive measure
at around $\t\approx 0.91$ and that the cyclic point of period three
appears at around $\t\approx 0.95$. The coordinates of the period
three, for $\t\approx 0.95$ point are approximately $(0.2,0.6,1)$.
In Section~\ref{locstabletheta1} we found that the period three
point exists as soon as inequality \eqref{percondit} is satisfied.
If $n=3$ and $\t=0.95$ then $g(1/\t-1)\approx 0.1900$,
$g^2(1/\t-1)\approx 0.5917$ and $1/2\t\approx 0.5263$. Hence, the
appearance of the period three point occurs a little before at the
threshold value of $\t$ predicted by inequality \eqref{percondit},
but it is of the required form $(g(1/\t-1),g^2(1/\t-1),1)$. This is
not surprising. We showed that a cyclic point of that form is stable
as soon as the inequality is satisfied. The eigenvalues vary
continuously with $\t$ so the point cannot suddenly become
unstable once $\t$ decreases below the threshold given in inequality
\eqref{percondit}.

\begin{figure}\label{feigdim45}
{\includegraphics[height=3in]{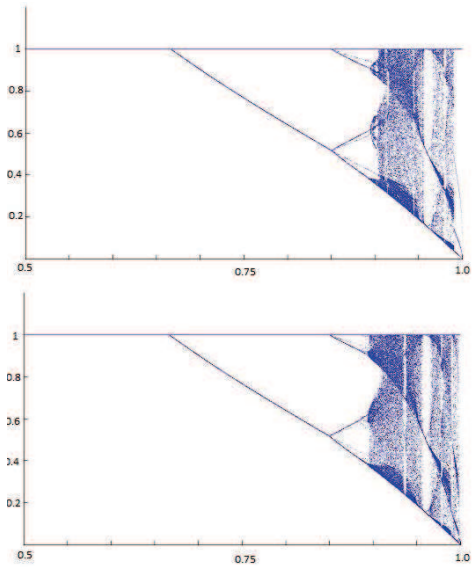}}
\caption{\label{bifudiag2} Feigenbaum diagram for the process on four
coordinates (top) and five coordinates (bottom).}
\end{figure}

The Feigenbaum diagrams for $n=4$ and $n=5$ are similar to the
diagram for $n=3$, and as it turns out that this holds in general
for all $n>3$. The main difference between $n=3$ and $n>3$ is the
appearance of a chaotic region for $0.95<\t<1$. It is remarkable that
a stable point of period three reappears around $\t\approx 0.95$.
For $n=4$ the stable cyclic point of period four appears at $\t
\approx 0.99$ and is still visible in this figure. For $n=5$ it
appears only at $\t\approx 0.999$ and it is not visible in this
picture. To show that our analysis holds
and that the periodic point does exist, we zoom in on step
sizes in $(0.95,1)$ in the next figure.

\begin{figure}\label{feigdimzoom}
{\includegraphics[height=1.5in]{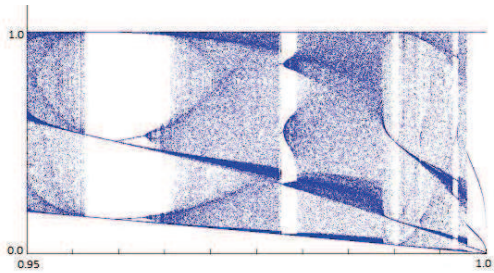}}
\caption{\label{bifudiag3} Feigenbaum diagram for the process on five
coordinates: zoom in on $\t\geq 0.95$.}
\end{figure}

The diagram shows the cyclic period five for
$\t$ near $1$.
This concludes our analysis of the Dikin process $w^k$. Now to
prove that this analysis makes sense, we still need to check that the
primal-dual AFS method displays the same type of chaotic behavior as
$w^k$. We will do that in the next and final section.

\section{Comparison to primal-dual AFS}

The iterative process $w^k$ has been derived by a linearization
of the primal-dual AFS method. To show that our bifurcation analysis
bears any relevance, we need to verify that a simlar route to chaor occurs in
actual LO problems. There is one complication.
The Dikin process involves a parameter $\t$ that
defines the step-size with respect to the maximum
$\alpha_{max}=\frac{\norm{xs}}{\max{xs}}$. So if we consider the
primal-dual AFS method, then we should set our step size
accordingly. This means that $\alpha$ should not be constant, which it
is in the \cite{}original primal-dual AFS method, but we should take it to be
equal to $\t\a_{\max}$. We modify the AFS method in this way and we put
$\a=\t\a_{max}$.

We take the same example as considered by Castillo and Barnes in
\cite{Barnes}:

\beq \begin{matrix} \min {10x_1+10x_2+5x_3+x_4-x_5}\\
\text{under the constraints}\\
x_1+2x_2-3x_3-2x_4-x_5=0\\
-x_1+2x_2-x_3-x_4-x_5=0\\
x\geq 0
\end{matrix}
\eeq

We take the same initial vectors $x_0$ and $y_0$ as Castillo and
Barnes and run \textit{our} modified primal-dual AFS method that we describe
in pseudo-code below. The numerical task of computing the limit of
the AFS process is not trivial, especially for a larger values of
the step size, because $x^k$ rapidly converges to zero which leads
to numerical problems, caused by inverting matrices that are ill conditioned.
Castillo and Barnes developed analytic formulas that
enabled them to still compute Feigenbaum diagrams with high precision.
Such an analytic exercise is beyond the scope
of our paper. We stop the computation once the duality gap reaches $10^{-10}$.

\begin{figure}[!hbt]
\begin{minipage}{4cm}
\hrule\par\bigskip
\centerline{{\Large \bf Modified Primal--Dual AFS}}
\par\bigskip
\par\hrule\par\bigskip

\parskip = 0pt
\par

\hspace{0.5cm}
\begin{minipage}[t]{4cm}
{\bf Parameters}\par \hspace{.2cm}
\begin{minipage}[t]{4cm}
$\varepsilon$ is the accuracy parameter;\\
$\t$ is the scaled step size;
\end{minipage}

\bigskip
{\bf Input}\par \hspace{.2cm}
\begin{minipage}[t]{4cm}
$(x^{0}, s^{0})$: the initial pair  of interior   feasible
solutions;
\end{minipage}

\bigskip
{\bf begin}\par \smallskip
\begin{minipage}[t]{4cm}
$x:=x^0;\, s:=s^0$; \\
\par\smallskip
\begin{minipage}[t]{4cm}
{\bf while}\quad $x^Ts>\varepsilon$ \quad {\bf do}
$w=xs$;\\
$\a_{max}=\frac{\norm{w}}{\max w}$;\\
$\a=\t\a_{max}$;\\
$x := x + \a\Dx$;\\
$y := y + \a\Dy$;\\
$s := s + \a\Ds$;
\end{minipage}
\par
\end{minipage}
\par \medskip
{\bf end.}
\end{minipage}
\par\bigskip
\par\hrule\par\bigskip
\end{minipage}
\caption{Primal--dual affine scaling algorithm with modified
step size $\alpha$. In our computations we put
$\varepsilon=10^{-10}$ and we plot results as soon as the duality
gap reaches values $\leq 0.001$}\label{fig:6}
\end{figure}
\begin{figure}[!hbt]
{\includegraphics[height=3in]{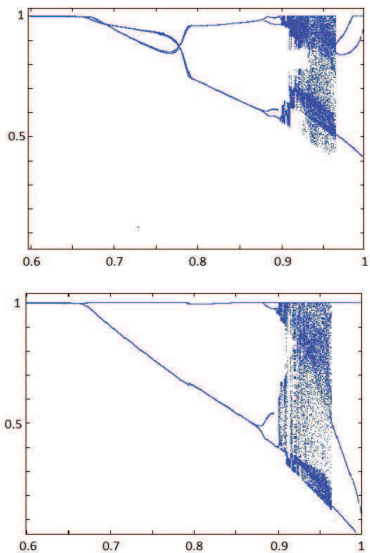}}
\caption{Feigenbaum diagrams for the Castillo-Barnes LO problem. Horizontal coordinate represents $\theta$.
Vertical axis contains the $\omega$-limit of
a coordinate of the scaled vector $w$. Second coordinate above. Fourth coordinate below.}\label{feigprimdual}
\end{figure}

We have computed the Feigenbaum diagram for the scaled process
$\frac {w^k}{\max w^k}$ that is given in Figure~\ref{feigprimdual}. The diagram below depicts the limit
of the fourth coordinate. There is a bifurcation for $\t=2/3$ and
another bifurcation close to $\t=0.86$, followed by a chaotic regime.
At the end of the diagram, for values of $\t$ close to $1$, we find a stable
periodic point.  This is similar to the diagrams that we computed
earlier for our process $w^k$, although the periodic point at the
end of the diagram is period three instead of period five. The
Feigenbaum diagram above, which depicts the second coordinate,
shows a different picture. The diagram bifurcates at $\t=2/3$ but the two
branches of the graph intersect twice between $2/3$ and $0.86$: once at
$\t\approx 0.69$ and once at $\t\approx 0.78$. At these values of $\t$,
the limit lands exactly on the unstable fixed point. We already noticed that this point is weakly repelling,
which is why the second coordinate has not yet fully converged to its $\omega$-limit
yet, even when the duality gap is $10^{-10}$.

The dual problem is degenerate
\beq \begin{matrix} \max {0}\\
\text{under the constraints}\\
y_1-y_2\leq 10\\
2y_1+2y_2\leq 10\\
-3y_1-y_2\leq 5\\
-2y_1-y_2\leq 1\\
-y_1-y_2\leq -1
\end{matrix}
\eeq

\begin{figure}[!hbt]
\resizebox{3in}{!}{\includegraphics[width=4in]{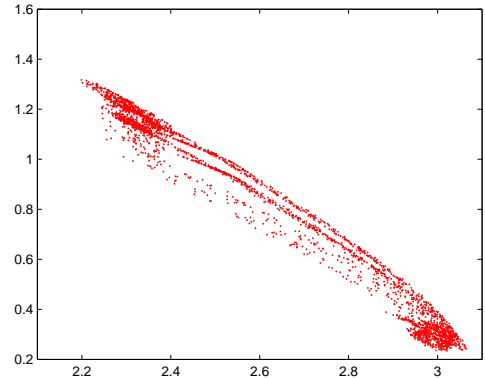}}
\caption{The omega-limit set of the vector $y$ in the dual problem for $\theta=0.94$ forms
a strange attractor in the feasible set.}\label{ylimit}
\end{figure}

All feasible points solve the dual problem. If $\t\leq 2/3$ then the
process $y^k$ converges to $(3.0513, 0.5522)$ but if $\t$ increases
beyond $2/3$ then the process no longer converges to a single point.
However, $y^k$ remains within the feasible set even for large values
of $\t$. Figure \ref{ylimit} contains the limit set that we
computed for $\t=0.94$. It has the contours of a H\'enon-like strange attractor.
The image of the attractor is slightly blurred since the orbit has not
fully converged yet.

It seems that the process $w^k$ that we have considered in this
paper represents the iterations of primal-dual AFS rather well.
We have tested other LO problems as well and we find similar Feigenbaum diagrams
for the vector $w$, regardless whether the dual problem is degenerate or not.
The algorithm converges to an optimal solution for
relatively high values of $\t$, so for a step-size that is close to
$\a_{max}$. This may indicate that a step-size that is larger than
$1/(15\sqrt n)$ is possible, if $\a$ is not taken to be constant but
is allowed to vary with $xs$, as in our computations.

\section{Conclusion} We have presented the Dikin process as an
archetype for general interior point methods.
The Dikin process is a one-parameter family with a route to chaos
that bears similarity to the
the logistic family, and which agrees with the chaotic
behaviour of interior point methods that has been previously
observed.

\end{document}